# An Asymptotic Variational Problem Modeling a Thin Elastic Sheet on a Liquid, Lifted at One End

David Padilla-Garza*

April 24, 2019

**Abstract**

We discuss a 1D variational problem modeling an elastic sheet on water, lifted at one end. Its terms include the membrane and bending energy of the sheet as well as terms due to gravity and surface tension. By studying a suitable Gamma-limit, we identify a parameter regime in which the sheet is inextensible, and the bending energy and weight of the sheet are negligible. In this regime, the problem simplifies to one with a simple and explicit solution.

## 1 Setting and model

This article models and analyzes an experiment in which a thin sheet on water is lifted at one end. Against expectations, the profile of the thin sheet on one side of the contact point, and the profile of the liquid gas interface on the other side of the contact point are exactly symmetric (see https://blogs.umass.edu/soft-matter/lecture-notes/ specifically the course by Benny Davidovitch, lecture 4 for experimental results, or (Deepak Kumar, 2019)). This is counter intuitive since the forces acting on the liquid gas interface are the gravitational pull of the liquid, and the effect of surface tension. On the other hand, the forces acting on the thin film are elastic forces, surface tension, and the gravitational pull of the liquid. Furthermore, the surface tension coefficients of the three different interfaces (liquid-gas, gas-solid, liquid-solid) are in principle different. This article provides a mathematical treatment of the problem: starting from well established first principles we deduce the solution, and prove that the profiles are symmetric.

The interaction of thin sheets with surface tension have been the focus of many recent works including (Neukirch, Antkowiak, & Marigo, 2013), and works dealing with wrinkling phenomena: (King, Schroll, Davidovitch, & Menon, 2012), (Huang, Davidovitch, Santangelo, Russell, & Menon, 2010), (Vella, Adda-Bedia, & Cerda, 2010), (Vella, Aussillous, & Mahadevan, 2004), (Paulsen et al., 2017). More broadly, (Vella & Mahadevan, 2005) and (Vella, Lee, & Kim, 2006) treat the interaction of surface tension with mechanical properties of materials at a millimeter scale.

*David Padilla-Garza. Courant Institute of Mathematical Sciences New York University 251 Mercer Street New York, N.Y. 10012. email: padilla@cims.nyu.edu. This work was partially supported by the National Science Foundation through grant DMS-1311833.

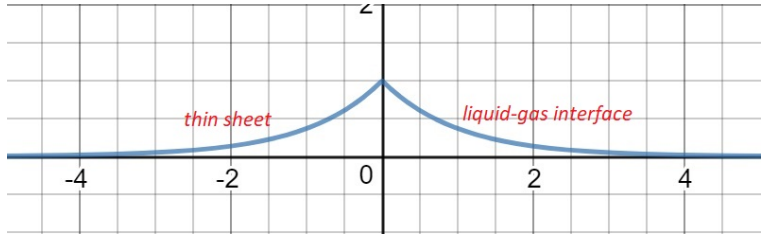

Figure 1: Schematic of the experimental observation showing the symmetry of the liquid-gas interface and thin sheet profiles. For the experimental results, see https://blogs.umass.edu/softmatter/lecture-notes/ specifically the course by Benny Davidovitch, lecture 4

To explain our setup, consider a thin elastic sheet of dimensions $[0, l] \times [0, l_1] \times [-\frac{h}{2}, \frac{h}{2}]$ suspended above the surface of a liquid. From now on, we identify the shape of the sheet with that of its midplane; $[0, l] \times [0, l_1] \times \{0\}$. Suppose that the right end $\{l\} \times [0, l_1] \times \{0\}$ is fixed at $\{\overline{x}, [0, l_1], \overline{y}\}$. Neglecting the effect of surface tension on the lateral edges, the problem has translation invariance and so we can describe the shape by a transversal cut to the midplane; let $(x(s), y(s))$ be the position of the point on the sheet $s \in [0, l]$. We choose a reference frame such that $y = 0$ corresponds to the average liquid level, i.e to a height such that no gravitational force due to liquid weight is exerted on the interface. We consider the problem where no external force is exerted on the point $(x(0), 0)$, which leads to a Neumann problem. Gravity will pull the sheet downwards and at the same time, surface tension will drive the surface of the liquid upwards. Our aim is to determine the shape of the sheet and the liquid-gas interface. Throughout the paper we will assume that $\gamma_{LG} > \gamma_{LS} + \gamma_{SG}$, that $\gamma_{LS} < \gamma_{LG} + \gamma_{GS}$ and that $\gamma_{GS} < \gamma_{LG} + \gamma_{LS}$ so that the sheet is in tension, there is wetting, and the model does not predict an ultra thin layer of liquid at the side of the sheet. We also assume that the contact set between the sheet and the liquid is an interval $(0, l)$.

We will use a variational model and hence the first step is to identify the physically realistic energy. It is clear that the energies involved in the setup are

- due to sheet weight
- due to liquid weight
- surface energy

- elastic energy (stretching+bending)

The energy due to the weight of the liquid and to its and surface tension is

$$E_1(y, y') = \int \frac{1}{2}\rho_L g y^2 + \gamma\sqrt{1 + (y')^2} dx \tag{1}$$

This choice is consistent with the well known Laplace Young equation, which describes the surface of a liquid $y(x)$ deformed by gravity and surface tension. In $1D$ it takes the form

$$\rho_L g y(x) = \gamma \frac{y''(x)}{(1 + (y'(x))^2)^{\frac{3}{2}}} \tag{2}$$

where

- $\rho_L$ is the density of the liquid
- $g$ is the gravitational acceleration
- $\gamma$ is the surface tension coefficient.

Equation 2 is the Euler-Lagrange equation for 1.

In what follows we will work with parametrized curves, and hence we must rewrite (1) for a parametrized curve $(x(s), y(s))$. This amounts to a change of variable, the new functional is

$$E_2(x, y, \dot{x}, \dot{y}) = \int \frac{1}{2}\rho_L g y^2 \dot{x} + \gamma\sqrt{(\dot{x})^2 + (\dot{y})^2} ds.$$

The energy of the sheet weight term is given by

$$h \int \rho_S g y(s) ds. \tag{3}$$

This choice is consistent with the well known catenary equation. This describes the shape of a string hanging by the effect of gravity with both endpoints fixed. Its solutions are critical points of the functional

$$E_3(y, y') = h \int \rho_S g y \sqrt{1 + (y')^2} dx,$$

where $\rho_S$ is the density of the sheet. This functional assumes that the string is inextensible. For an extensible string of reference length $l$ and shape given by $(x(s), y(s)), \ s \in [0, l]$ the energy is given by

$$E_3(y, y') = h \int \frac{1}{\nu}\rho_S g y \sqrt{1 + (y')^2} dx,$$

where $\nu = \sqrt{\dot{x}^2 + \dot{y}^2}$. Therefore for a parametrized curve $(x(s), y(s))$ the energy is given by (3).

Lastly, the elastic energy of our thin elastic sheet $(x(s), y(s))$ can be approximated up to order $h^3$ by

$$E_4 = Eh \int (\sqrt{\dot{x}^2 + \dot{y}^2} - 1)^2 + ch^2\kappa^2 ds,$$

where $h$ is the thickness of the sheet, $E$ is the elastic modulus and $\kappa$ is the curvature, and $c$ is a suitable constant, which we will take henceforth to be 1 for ease of notation.

We are now ready to write the energy functional. Divide the setup into four vector valued functions $(w_1(s), z_1(s)), (x_1(s), y_1(s)), (w_2(s), z_2(s)), (x_1(s), y_1(s))$, where

- $(w_1(s), z_1(s))$ describes the liquid-gas interface to the left of the sheet.

- $(x_1(s), y_1(s))$ describes the shape of the sheet to the left of the contact point.

- $(x_2(s), y_2(s))$ describes the shape of the sheet to the right of the contact point.

- $(w_2(s), z_2(s))$ describes the shape of the liquid-gas interface to the right of the contact point.

See figure 1. We assume that the sheet and the liquid-gas interface meet in a single point.

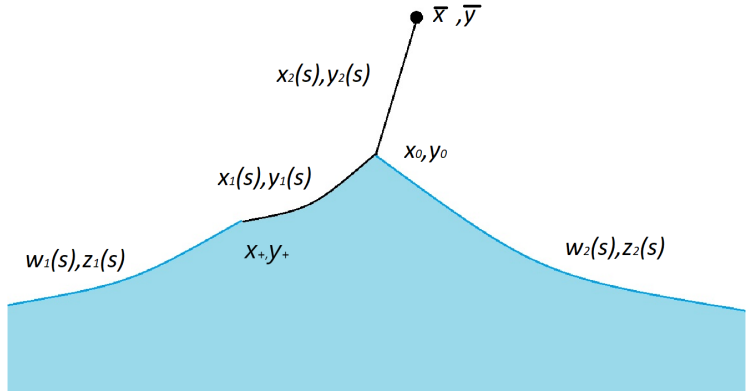

Figure 2: Rough sketch of of experiment.

The energy functional is

$$\min_{(x_+,y_+),(x_0,y_0),l} \min_A \int_0^\infty \gamma_{LG}\sqrt{\dot{z}_1^2+\dot{w}_1^2} + \rho_L g \frac{1}{2}z_1^2|\dot{w}_1|ds$$

$$+ \int_0^l (\gamma_{SG}+\gamma_{SL})\sqrt{\dot{x}_1^2+\dot{y}_1^2} + \rho_L g \frac{1}{2}y_1^2\dot{x}_1 + gh\rho_S y_1$$

$$+ Eh(\sqrt{\dot{x}_1^2+\dot{y}_1^2}-1)^2 + Eh^3\kappa^2 ds$$

$$+ \int_0^{L-l} 2\gamma_{SG}\sqrt{\dot{x}_2^2+\dot{y}_2^2} + gh\rho_S y_2$$

$$+ Eh(\sqrt{\dot{x}_2^2+\dot{y}_2^2}-1)^2 + Eh^3\kappa^2 ds$$

$$+ \int_0^\infty \gamma_{LG}\sqrt{\dot{z}_2^2+\dot{w}_2^2} + \rho_L g \frac{1}{2}z_2^2|\dot{w}_2|ds$$

where

$$\begin{aligned}A = \{&(w_1,z_1),(x_1,y_1),(w_2,z_2),(x_2,y_2)|\\
&(w_1(0),z_1(0))=(x_1(0),y_1(0))=(x_+,y_+)\\
&(x_2(0),y_2(0))=(x_1(l),y_1(l))=(w_2(0),z_2(0))=(x_0,y_0) \quad (4)\\
&(x_2(L-l),y_2(L-l))=(\overline{x},\overline{y})\\
&(\dot{x}_2(0),\dot{y}_2(0))=(\dot{x}_1(l),\dot{y}_1(l))\}\end{aligned}$$

(note that the last condition is necessary to avoid concentration of bending energy at the point $(x_1(l),y_1(1)))$. The next step is to write the model in dimensionless coordinates. Let $\tilde{s}$ be the original position in reference coordinates and $(\tilde{x}_i,\tilde{y}_i),(\tilde{w}_i,\tilde{z}_i)$ be a setup as described above. Let $s=\frac{\tilde{s}}{L}$, $(x_i,y_i)=\frac{1}{L}(\tilde{x}_i,\tilde{y}_i)$, $(w_i,z_i)=\frac{1}{L}(\tilde{w}_i,\tilde{z}_i)$ and $\hat{h}=\frac{h}{L}$. Then (taking $L=1$ in the definition of $A$) the functional in the dimensionless variables is

$$\min_{(x_+,y_+),(x_0,y_0),l} \min_A \int_0^\infty \gamma_{LG}L\sqrt{\dot{z}_1^2+\dot{w}_1^2} + \rho_L gL^3\frac{1}{2}z_1^2|\dot{w}_1|ds$$

$$+ \int_0^l L(\gamma_{SG}+\gamma_{SL})\sqrt{\dot{x}_1^2+\dot{y}_1^2} + \rho_L gL^3\frac{1}{2}y_1^2\dot{x}_1 + \rho_S gL^3\hat{h}y_1$$

$$+ EL^2\hat{h}(\sqrt{\dot{x}_1^2+\dot{y}_1^2}-1)^2 + E\hat{h}^3L^2\kappa^2 ds$$

$$+ \int_0^{1-l} 2\gamma_{SG}L\sqrt{\dot{x}_2^2+\dot{y}_2^2} + \rho_S gL^3\hat{h}y_2 + E\hat{h}L^2(\sqrt{\dot{x}_2^2+\dot{y}_2^2}-1)^2$$

$$+ E\hat{h}^3L\kappa^2 ds$$

$$+ \int_0^\infty (\gamma_{LG})L\sqrt{\dot{z}_2^2+\dot{w}_2^2} + \rho_L gL^3\hat{h}\frac{1}{2}z_2^2|\dot{w}_2|ds.$$

Normalizing by $\hat{h}EL^2$ we obtain

$$\min_{(x_+,y_+),(x_0,y_0),l} \min_A \int_0^\infty \frac{\gamma_{LG}}{E\hat{h}L}\sqrt{\dot{z}_1^2 + \dot{w}_1^2} + \frac{\rho_L g L}{E\hat{h}}\frac{1}{2}z_1^2|\dot{w}_1|ds$$
$$+ \int_0^l \frac{(\gamma_{SG} + \gamma_{SL})}{E\hat{h}L}\sqrt{\dot{x}_1^2 + \dot{y}_1^2} + \frac{\rho_L g L}{E\hat{h}}\frac{1}{2}y_1^2\dot{x}_1 + \frac{\rho_S g L}{E}y_1$$
$$+ (\sqrt{\dot{x}_1^2 + \dot{y}_1^2} - 1)^2 + \hat{h}^2\kappa^2 ds$$
$$+ \int_0^{1-l} \frac{2\gamma_{SG}}{E\hat{h}L}\sqrt{\dot{x}_2^2 + \dot{y}_2^2} + \frac{\rho_L g L}{E}\hat{h}y_2$$
$$+ (\sqrt{\dot{x}_2^2 + \dot{y}_2^2} - 1)^2 + \hat{h}^2\kappa^2 ds$$
$$+ \int_0^\infty \frac{\gamma_{LG}}{E\hat{h}L}\sqrt{\dot{z}_2^2 + \dot{w}_2^2} + \frac{\rho_L g L}{E\hat{h}}\frac{1}{2}z_2^2|\dot{w}_2|ds.$$

We introduce the new constants

- $\frac{\gamma_i}{E\hat{h}L} = A_i$ with $i = LG$, $i = SL$ or $i = SG$
- $\frac{g\rho_S L}{E} = B$
- $\frac{g\rho_L L}{E\hat{h}} = C$.

This functional is identically $+\infty$ due to the first and fourth term. In order to get a well defined functional, we subtract from these terms the energy of the state of zero energy, which we define to be a horizontal semi-infinite liquid-gas interface with height 0 starting at $\overline{x}$. The resulting functional becomes:

$$\min_{(x_+,y_+),(x_0,y_0),l} \min_A \int_0^\infty A_{LG}(\sqrt{\dot{z}_1^2 + \dot{w}_1^2} + \dot{w}_1) + C\frac{1}{2}z_1^2|\dot{w}_1|ds + A_{LG}w_1(0)$$
$$+ \int_0^l (A_{SG} + A_{SL})\sqrt{\dot{x}_1^2 + \dot{y}_1^2} + C\frac{1}{2}y_1^2\dot{x}_1 + By_1$$
$$+ (\sqrt{\dot{x}_1^2 + \dot{y}_1^2} - 1)^2 + \hat{h}^2\kappa^2 ds$$
$$+ \int_0^{1-l} 2A_{SG}\sqrt{\dot{x}_2^2 + \dot{y}_2^2} + By_2 + (\sqrt{\dot{x}_2^2 + \dot{y}_2^2} - 1)^2 + \hat{h}^2\kappa^2 ds$$
$$+ \int_0^\infty A_{LG}(\sqrt{\dot{z}_2^2 + \dot{w}_2^2} - \dot{w}_2) + C\frac{1}{2}z_2^2|\dot{w}_2|ds - A_{LG}w_2(0). \tag{5}$$

## 2 Analysis and solution

This functional looks awfully cumbersome. In order to obtain a more tractable functional, we introduce the function

$$\phi^\pm(x_0, y_0) = \min_{x(0)=y_0} \int_0^\infty E_{graph}^{LY}(y, y')dx \mp A_{LG}x_0 \tag{6}$$

where
$$E_{graph}^{LY}(y,y') = A_{LG}(\sqrt{1+(y'(x))^2} - 1) + C\frac{1}{2}y^2.$$

We also introduce the curve
$$(x(s), y(s)) = \begin{cases} (x_1(s), y_1(s)) & s \in [0, l] \\ (x_2(s-l), y_2(s-l)) & s \in [l, 1]. \end{cases}$$

Before we can proceed, we need some basic facts about $\phi$, which will be proved in the appendix.

**Lemma 1.** *There exists a $y^*$ (which depends on $A_{LG}$ and $C$) such that for every $y^* \geq y_0 > 0$ there exists a unique solution $y(x)$ to the Laplace Young equation such that $y(0) = y_0$ and $\lim_{x \to \infty} y(x) = 0$.*

The next step is to prove uniqueness of solutions among functions that are not necessarily graphs. This amounts to proving uniqueness of solutions to the Laplace-Young functional written for parametric curves.

**Proposition 2.** *The problem*
$$\min \int_0^\infty \frac{1}{2} C y^2 |\dot{x}| + A_{LG}(\sqrt{(\dot{x})^2 + (\dot{y})^2} - \dot{x}) ds \tag{7}$$

*among $x, y$ such that $(x(0), y(0))$ is fixed and $\lim_{s \to \infty}(x(s), y(s)) = (\infty, 0)$ has a solution that is unique up to reparametrization.*

**Proposition 3.** *The function $\phi^\pm$ (defined by (6)) is smooth away from $y*$.*

The functional becomes
$$E_{\hat{h}}(x, y, \dot{x}, \dot{y}, l) = \int_0^1 (\mathbf{1}_{[0,l]}(s)[A_{LS} + A_{SG}] + \mathbf{1}_{[l,1]}(s) 2 A_{SG}) \sqrt{\dot{x}_1^2 + \dot{y}_1^2}$$
$$+ C\mathbf{1}_{[0,l]}(s)\frac{1}{2} y_1^2 \dot{x}_1 + By_1 + (\sqrt{\dot{x}^2 + \dot{y}^2} - 1)^2 + \hat{h}^2 \kappa^2 ds$$
$$+ \phi^-(x(0), y(0)) + \phi^+(x(l), y(l)).$$

We now study the limit of the functional as the nondimensionalized thickness $\hat{h}$ tends to 0. We assume that the parameters have the scaling

- $\lim_{\hat{h} \to 0} \frac{1}{\hat{h}^\alpha} A_i^{\hat{h}} \to A_i^*$ with $i = LG$, $i = SG$ or $i = LS$
- $\lim_{\hat{h} \to 0} \frac{1}{\hat{h}^\alpha} C^{\hat{h}} \to C^*$
- $\lim_{\hat{h} \to 0} \frac{1}{\hat{h}^{\alpha+\epsilon}} B^{\hat{h}} \to B^*$,

for $\alpha \in (0, 2)$ and $\epsilon > 0$. While the scalars $A_i^{\hat{h}}, B^{\hat{h}}, C^{\hat{h}}$ depend on $\hat{h}$, this dependence will be henceforth suppressed for convenience of notation. From now on we drop the superscript on $\hat{h}$, so that $h$ stands for the nondimensionalized thickness.

**Theorem 4.** *Let $\mathcal{E}_h$ be the functional*

$$\mathcal{E}_h(x,y,\dot{x},\dot{y},l) = \begin{cases} \bullet \dfrac{1}{h^\alpha}\left(\int_0^1 \left(\mathbf{1}_{[0,l]}(s)[A_{LS}+A_{SG}] + \mathbf{1}_{[l,1]}(s)2A_{SG}\right)\sqrt{\dot{x}^2+\dot{y}^2}\right. \\ \qquad\qquad \left. + C\mathbf{1}_{[0,l]}(s)\dfrac{1}{2}y^2\dot{x}^2 + By + (\sqrt{\dot{x}+\dot{y}^2}-1)^2 + h^2\kappa^2(s)ds\right) \\ \qquad + \phi^-(x(0),y(0)) + \phi^+(x(l),y(l)) \\ \qquad \text{if } (x(s),y(s)) \in W^{2,2}[0,1] \\ \bullet \infty \qquad \text{otherwise.} \end{cases}$$
(8)

*then $\mathcal{E}_h$ $\Gamma$ converges (in the weak\* $W^{1,\infty}$ topology) to the functional*

$$E(x,y,\dot{x},\dot{y},l) = \begin{cases} \bullet \int_0^1 (\mathbf{1}_{[0,l]}(s)[A^*_{LS}+A^*_{SG}] + \mathbf{1}_{[l,1]}(s)2A^*_{SG}) \\ \qquad + C^*\mathbf{1}_{[0,l]}(s)\dfrac{1}{2}y^2\dot{x}ds + \phi^-(x(0),y(0)) + \phi^+(x(l),y(l)) \\ \qquad \text{if } (x,y) \in W^{1,\infty}[0,1], \dot{x}^2+\dot{y}^2 \leq 1 \\ \bullet \infty \qquad \text{otherwise.} \end{cases}$$
(9)

Before embarking on the proof, let us make a few remarks about the differences between these functionals:

- The functionals $\mathcal{E}_h$ charge infinite bending energy to corners, while the limit $E$ charges zero energy to corners.

- The functionals $\mathcal{E}_h$ charge a positive but finite amount of membrane energy to compression and stretching, while the limit $E$ charges zero membrane energy to compression and infinite energy to stretching.

- In the functionals $\mathcal{E}_h$ the surface energy is proportional to deformed arclength, while in the limit $E$ the surface energy is proportional to reference arclength. The physical interpretation is that in an ultra thin sheet, a contracted configuration consists of fine wrinkles, hence the surface energy is proportional to the area of the wrinkled sheet.

We will use two lemmas in the proof, the first is quite standard, while the proof of lemma 7 can be found in the appendix.

**Lemma 5.** *If $f_n, f \in W^{1,\infty}[0,1]$ and $f_n \to f$ pointwise with $\|f'_n\|_{L^\infty} \leq K < \infty$, $\|f'\|_{L^\infty} \leq K < \infty$ then $f_n \to f$ weak\* in $W^{1,\infty}[0,1]$.*

The next lemma is a very intuitive and elementary fact about isometric embeddings. Similar results can be found for example in (Bella & Kohn, 2014) and (Conti & Maggi, 2008).

**Definition 6.** *A function $f:(a,b) \to \mathbf{R}^n$ is a* short embedding *if $\|\dot{f}\| \leq 1$.*

**Lemma 7.** *Let $f : [0,1] \to \mathbf{R}^2$ be a $C^\infty$ short embedding. For every $\epsilon > 0$ there exists $g \in C^\infty[0,1]$ such that*

- *The function $g$ is an isometric embedding.*
- *$|f(s) - g(s)| < \epsilon$ for all $s \in [0,1]$.*

With these lemmas, we now give the proof of theorem 4

*Proof.* (Of Theorem 4 We start by proving the lim inf inequality. Let $(x_n, y_n) \in W^{1,\infty}[0,1]$ converge weak* to $(x, y) \in W^{1,\infty}$ with $\dot{x}^2 + \dot{y}^2 \leq 1$ a.e. and let $l_n \to l$.

We get immediately that

$$\int_0^1 (\mathbf{1}_{[0,l_n]}(s)[A_{LS} + A_{SG}] + \mathbf{1}_{[l_n,1]}(s) 2 A_{SG}^*) ds$$

converges to

$$\int_0^1 (\mathbf{1}_{[0,l]}(s)[A_{LS} + A_{SG}] + \mathbf{1}_{[l,1]}(s) 2 A_{SG}^*) ds. \tag{10}$$

Because of the definition of weak* $L^\infty$ convergence we have

$$|(x_n(s), y_n(s)) - (x(s), y(s))| = |\int [(\dot{x}_n, \dot{y}_n) - (\dot{x}, \dot{y})] \mathbf{1}_{(s,1)}|$$
$$\to 0$$

for every $s \in [0,1]$, in particular and since $\phi^\pm$ is continuous, we have

$$\phi^-(x_n(0), y_n(0)) \to \phi^-(x(0), y(0)).$$

We also have

$$|\phi^+(x_n(l_n), y_n(l_n)) - \phi^+(x(l), y(l))| \leq |\phi^+(x_n(l_n), y_n(l_n)) - \phi^+(x_n(l), y_n(l)| + \\ |\phi^+(x_n(l), y_n(l)) - \phi^+(x(l), y(l)|$$

and

$$|(x_n(l_n), y_n(l_n)) - (x_n(l), y_n(l))| \leq |l - l_n|,$$

hence $x_n(l_n), y_n(l_n) \to (x(l), y(l)$ and $\phi^+(x_n(l_n), y_n(l_n)) \to \phi^+((x(l), y(l))$.

Lastly, by the weak-strong lemma, $y_n^2 \dot{x}_n \to y^2 \dot{x}$ weak* in $L^\infty$, and hence

$$\int y_n^2 \dot{x}_n \mathbf{1}_{(0,l)} \to \int y^2 \dot{x} \mathbf{1}_{(0,l)}.$$

We also have that $\| \int y_n^2 \dot{x}_n \|_{L^\infty} \leq K$, hence by Hölder we have

$$\left| \int y_n^2 \dot{x}_n \mathbf{1}_{(0,l_n)} - \int y^2 \dot{x} \mathbf{1}_{(0,l)} \right| \leq \left| \int (y_n^2 \dot{x}_n - y^2 \dot{x}) \mathbf{1}_{(0,l)} \right| +$$

$$\left| \int (y_n^2 \dot{x}_n - y^2 \dot{x})(\mathbf{1}_{(0,l)} - \mathbf{1}_{(0,l)}) \right|$$

$$\leq \left| \int (y_n^2 \dot{x}_n - y^2 \dot{x}) \mathbf{1}_{(0,l)} \right| + K|l - l_n|$$

$$\to 0.$$

For the other terms, we assume that $\liminf \mathcal{E}_{h_n}(x_n, y_n, l_n) < \infty$ since otherwise inequality is trivial. Then we have

$$\int_a^b (\sqrt{\dot{x}_n^2 + \dot{y}_n^2} - 1)^2 ds \to 0$$

for every $(a, b)$. Hence by Hölder's inequality we have

$$\int_a^b \left( \sqrt{\dot{x}_n^2 + \dot{y}_n^2} - 1 \right) ds \leq \sqrt{(b-a) \int_a^b (\sqrt{\dot{x}_n^2 + \dot{y}_n^2} - 1)^2 ds}$$

$$\to 0$$

hence

$$\int_a^b \sqrt{\dot{x}_n^2 + \dot{y}_n^2} \, ds \to b - a.$$

Using Hölder once again we get

$$\int \mathbf{1}_{(l_n,l)} (\sqrt{\dot{x}_n^2 + \dot{y}_n^2} - 1) \leq \sqrt{|l_n - l|} \, \| \sqrt{\dot{x}_n^2 + \dot{y}_n^2} - 1 \|_{L^2[0,1]}$$

$$\to 0.$$

By the triangle inequality,

$$\int \mathbf{1}_{(l_n,l)} \sqrt{\dot{x}_n^2 + \dot{y}_n^2} \leq \left| \int \mathbf{1}_{(l_n,l)} (\sqrt{\dot{x}_n^2 + \dot{y}_n^2} - 1) \right| + |l_n - l|$$

$$\to 0.$$

We thus get

$$\int \mathbf{1}_{(0,l_n)} \sqrt{\dot{x}_n^2 + \dot{y}_n^2} ds = \int \mathbf{1}_{(0,l)} \sqrt{\dot{x}_n^2 + \dot{y}_n^2} ds + \int \mathbf{1}_{(l,l_n)} \sqrt{\dot{x}_n^2 + \dot{y}_n^2} ds$$

$$\to l.$$

Note that since $\dot{x}^2 + \dot{y}^2$ is l.s.c with respect to weak* $W^{1,\infty}$ convergence we get that $\dot{x}^2 + \dot{y}^2 \leq 1$ a.e. This completes the proof of the lim inf inequality.

We now turn to the lim sup inequality. Let $(x, y) \in W^{1,\infty}[0, 1]$ with $\dot{x}^2 + \dot{y}^2 \leq 1$ a.e. Let $x_n = x * \mu_{\frac{1}{n}}$ and $y_n = y * \mu_{\frac{1}{n}}$ where $\mu : \mathbf{R} \to \mathbf{R}^+$ is an even $C^\infty$ function such that $\text{supp}[\mu] \in [-1, 1]$ and $\int \mu(x)dx = 1$, and $\mu_{\frac{1}{n}} = n\mu(nx)$. In order to have $(x_n, y_n)$ well defined in the intervals $(0, \frac{1}{n})$ and $(1 - \frac{1}{n}, 1)$ we can extend them to functions on the interval $(-1, 2)$ by odd reflexion around 0 and 1. This way $(x_n(0), y_n(0)) = (x(0), y(0))$ and $(x_n(1), y_n(1)) = (x(0), y(0))$. Since $(\dot{x}, \dot{y}) \in B(0, 1)$ and $B(0, 1)$ is a convex set, we have that $(\dot{x}_n, \dot{y}_n) \in B(0, 1)$ and $(x_n(s), y_n(s))$ is a short embedding. By classical harmonic analysis, $(x_n, y_n) \to (x, y)$ in $C^0[0, 1]$, and since $\dot{x}, \dot{x}_n, \dot{y}, \dot{y}_n \leq 1$ we have that $(x_n, y_n) \to (x, y)$ weak* in $W^{1,\infty}$. We also have that $(x_n, y_n) \in C^\infty$.

Using Lemma 7, there exists a $C^\infty$ isometric immersion $\hat{x}_n, \hat{y}_n : [0, 1] \to \mathbf{R}^2$ such that $\| (\hat{x}_n, \hat{y}_n) - (x_n, y_n) \|_{C^0} < \frac{1}{n}$. Since $\max\{\dot{\hat{x}}_n, \dot{\hat{y}}_n, \dot{x}_n, \dot{y}_n\} \leq 1$ we have that $(\hat{x}_n, \hat{y}_n) - (x_n, y_n) \to 0$ weak* in $W^{1,\infty}$.

Define a sequence $\sigma(n)$ as
$$\sigma(1) = 1$$
and
$$\sigma(n+1) = \begin{cases} \sigma(n) + 1 & \text{if } h_n^{2-\alpha} \int \kappa^2(\hat{x}_{\sigma(n)+1}, \hat{y}_{\sigma(n)+1}, s)ds \leq \frac{1}{\sigma(n)+1} \\ \sigma(n) & \text{if not} \end{cases}$$

where $\kappa(\hat{x}_n, \hat{y}_n, s)$ is the Gaussian curvature of the curve $(\hat{x}_n, \hat{y}_n)$ at point $s$.

The approximating sequence $(\theta_n(s), \gamma_n(s))$ is
$$(\theta_n(s), \gamma_n(s)) = (\hat{x}_{\sigma(n)}, \hat{y}_{\sigma(n)}).$$

Note that $\sigma(n) \to \infty$ since $h_n^{2-\alpha} \int \kappa^2(\hat{x}_m, \hat{y}_m, s)ds \to 0$ as $n \to \infty$ for a fixed $m$. Hence $(\theta_n(s), \gamma_n(s)) \to (x(s), y(s))$ in weak* $W^{1,\infty}$. Therefore (arguing as in the liminf inequality) we have that

$$\frac{1}{h_n^\alpha} \left( \int_0^1 \left( \mathbf{1}_{[0,l]}(s)[A_{LS} + A_{SG}] + \mathbf{1}_{[l,1]}(s)2A_{SG} \right) \sqrt{\dot{\theta}_n^2 + \dot{\gamma}_n^2} \, ds \right)$$
$$+ \frac{1}{h_n^\alpha} \left( C\mathbf{1}_{[0,l]}(s)\frac{1}{2}(\gamma_n)^2\dot{\theta}_n + B\gamma_n ds \right)$$
$$+ \phi^-(\theta_n(0), \gamma_n(0)) + \phi^+(\theta_n(l), \gamma_n(l))$$

converges to (since $\sqrt{\dot{\theta}_n^2 + \dot{\gamma}_n^2} = 1$)

$$\frac{1}{h_n^\alpha} \left( \int_0^1 \left( \mathbf{1}_{[0,l]}(s)[A_{LS}^* + A_{SG}^*] + \mathbf{1}_{[l,1]}(s)2A_{SG}^* \right) + C^*\mathbf{1}_{[0,l]}(s)\frac{1}{2}(x)^2\dot{y} + B^*y ds \right)$$
$$+ \phi^-(x(0), y(0)) + \phi^+(x(l), y(l)).$$

We also have that
$$h_n^{2-\alpha} \int_0^1 \kappa^2(\gamma_n, \theta_n, s)ds \to 0;$$

therefore
$$\limsup \mathcal{E}_n(\theta_n, \gamma_n) \leq \mathcal{E}(x,y)$$

And the proof is complete. □

**Proposition 8.** *Problem* (9) *has at least one solution, and it is the limit (in the weak $W^{1,2}$ topology) of a subsequence of minimizers of the functional $\mathcal{E}_h$.*

**Remark 9.** *1. This result does not immediately follow from $\Gamma$ convergence, since that result used the weak\* $W^{1,\infty}$. However, the same ideas applied to the weak $W^{1,2}$ topology yield the result.*

*2. We will prove later that the $\Gamma$ limit has a unique minimizer, therefore we can use the fact that every subsequence of minimizers has a further subsequence which converges to the minimizer of E to get that the original sequence converges to the minimizer.*

*3. Proposition 8 also holds for almost minimizers.*

*Proof.* Take a subsequence $(x_n, y_n, l_n)$ such that
$$\mathcal{E}_n(x_n, y_n, l_n) \leq m_n + \frac{1}{n},$$
where
$$m_n = \inf_A \mathcal{E}_n,$$
then
$$(\sqrt{\dot{x}^2 + \dot{y}^2} - 1)^2 \leq C < \infty,$$
and therefore $(x_n, y_n, l_n)$ converge (up to a subsequence) to a limit $(x, y, l)$ in $W^{1,2} \times \mathbf{R}$. The proposition will be proved if we prove that
$$\mathcal{E}(x, y, l) \leq \liminf \mathcal{E}(x_n, y_n, l_n),$$
in order to do so, we just need to make use of standard results in Sobolev spaces:

- Let $f_n : \mathbf{R} \to \mathbf{R}$ measurable and uniformly bounded by $K < \infty$ a.e. If $f_n \to f$ in $L^1$ then $f_n \to f$ in $L^2$.

- If $f_n \to f$ in $L^2$ then $f_n^2 \to f^2$ in $L^1$.

- If $(x_n, y_n) \to (x, y)$ weakly in $W^{1,2}[0,1]$ and $x_n(0), y_n(0) = (0,0)$ for each $n$ then $y_n^2 \dot{x}_n \to y^2 \dot{x}$ weakly in $L^2$.

Note the other terms in the functional do not need an essentially new proof, since the same ideas in the $\Gamma$ lim inf inequality prove they are l.s.c. □

We continue with an elementary observation about the functional which will significantly simplify the study of the problem.

**Remark 10.** *(Invariance under reparametrization). Let $\tau : [0,1] \to [0,1]$ be a diffeomeorphism such that $\tau(0) = 0$, $\tau(1) = 1$ and $\tau(l) = l$. Let $z = \tau(s)$, then*

$$\int_0^1 (\mathbf{1}_{[0,l]}(s)[A_{LS}^* + A_{SG}] + \mathbf{1}_{[l,1]}(s)2A_{SG}^*) + C^*\mathbf{1}_{[0,l]}(s)\frac{1}{2}y_1^2 \frac{d}{ds}x_1 ds = \\ \int_0^1 (\mathbf{1}_{[0,l]}(z)[A_{LS}^* + A_{SG}] + \mathbf{1}_{[l,1]}(z)2A_{SG}^*) + C^*\mathbf{1}_{[0,l]}(z)\frac{1}{2}y_1^2 \frac{d}{dz}x_1 dz \tag{11}$$

*Proof.* Chain rule. □

**Proposition 11.** *The minimizer of $E$ has $\| (\dot{x}, \dot{y}) \| = 1$ and therefore the deformed sheet has length $1$.*

*Proof.* The idea of this proof is that, by stretching the sheet a little bit along $(w_1, z_1)$, we relieve a little bit of the surface energy of $(w_1, z_1)$ without adding any gravitational or surface energy. Take $(x(s), y(s), l)$ such that length$(x, y) = h < 1$. For simplicity, assume $(x(s), y(s))$ has constant speed $h = \sqrt{\dot{x}^2 + \dot{y}^2}$. Take $(w_1, z_1)$ to be the solution of problem (7) parametrized by arclength. For any $\tau \in (0, 1-h)$ we can define.

$$(x_\tau(s), y_\tau(s)) = \begin{cases} x(\frac{h+\tau}{h}s), y(\frac{h+\tau}{h}s) \text{ if } s \in (0, \frac{h}{h+\tau}) \\ w_1(s - \frac{h}{h+\tau}), z_1(s - \frac{h}{h+\tau}) \text{ if } s \in (\frac{h}{h+\tau}, 1); \end{cases}$$

with $l_\tau = (\frac{h}{h+\tau})l$. Recalling the definition of $\phi$ we have that $E(x, y, \dot{x}, \dot{y}, l) - E(x_\tau, y_\tau, \dot{x}_\tau, \dot{y}_\tau, l_\tau) = A_{LG}(1 - \frac{h}{h+\tau})$, by letting $\tau$ tend to 0, we obtain that $x, y$ is not a local minimizer. Hence, minimizer(s) are stretched.

□

It would be possible to state the problem as minimizing (9) among all $(x, y)$ such that $\| (\dot{x}, \dot{y}) \| = 1$. But with the help of the last proposition 11, it is possible to solve our problem by taking perturbations that fix a contact point and the arclength at which the contact point occurs. This results in two lagrange multipliers. The following proposition allows us to reduce that even further to a single lagrange multiplier. To do this we consider the functional $E$ with the constraint of pointwise isometry relaxed to a constraint on arclength.

**Proposition 12.** *Let $x, y, l$ be a minimizer of*

$$\mathcal{E}(x, y, l) = \int_0^1 (\mathbf{1}_{[0,l]}(s)[A_{LS}^* + A_{SG}^*] + \mathbf{1}_{[l,1]}(s)2A_{SG}^*) \\ + C^*\mathbf{1}_{[0,l]}(s)\frac{1}{2}y^2 \dot{x} ds + \phi^-(x(0), y(0)) + \phi^+(x(l), y(l)) \tag{12}$$

*among curves $(x(s), y(s))$ such that $\sqrt{\dot{x} + \dot{y}^2} = 1$ (the energy functional (12) is (9) without the pointwise constraint of short isometry), then it minimizes $\mathcal{E}$ with respect to all perturbations with length one, not necessarily having unit speed.*

*Proof.* Let $(x(s), y(s), l)$ be the minimizer, and $\tilde{x}, \tilde{y}, \tilde{l}$ any other curve such that $\tilde{l} \in [0,1]$ and length$\tilde{x}, \tilde{y} = 1$. Let $\varphi : [0,1] \to [0,1]$ be the reparamatrization by arclength of $(\tilde{x}, \tilde{y})$ and $l^* = \varphi(\tilde{l})$. Let $x^*(s), y^*(s) = \tilde{x}(\varphi(s), \tilde{y}(\varphi(s)))$. By (10)

$$\mathcal{E}(x^*, y^*, l^*) = \mathcal{E}(\tilde{x}, \tilde{y}, \tilde{l}),$$

also by definition of $(x, y, l)$ we have

$$\mathcal{E}(x^*, y^*, l^*) \geq \mathcal{E}(x, y, l)$$

□

We are finally in a position to state the problem as minimizing the energy functional with a single Lagrange multiplier. We can use perturbations that fix $(x(l), y(l))$ to deduce the Euler-Lagrange equations.

It is standard in the Calculus of Variations that the critical point(s) of (11) among functions fixed at 0 and $l$, with fixed arclength 1 satisfy the same equation as those of the functional

$$\widetilde{E}(x, y, \dot{x}, \dot{y}, l) = \int_0^1 (\mathbf{1}_{[0,l]}(s)[A^*_{LS} + A^*_{SG}] + \mathbf{1}_{[l,1]}(s) 2 A^*_{SG}) \sqrt{\dot{x}_1^2 + \dot{y}_1^2} +$$
$$C^* \mathbf{1}_{[0,l]}(s) \frac{1}{2} y_1^2 \dot{x}_1 + \lambda(\sqrt{\dot{x}_1^2 + \dot{y}_1^2} - 1) ds +$$
$$\phi^-(x(0), y(0)) + \phi^+(x(l), y(l))$$

for some $\lambda \in \mathbf{R}$ (the Lagrange multiplier), among test functions that fix $(x(0), y(0))$ and $(x(l), y(l))$.

**Proposition 13.** *Critical points of the functional* (5) *satisfy the weak form of the Euler Lagrange equation*

$$\frac{d}{ds}\left[(\lambda + A_{LS} + A_{SG})\frac{\dot{x}_1}{\sqrt{\dot{x}_1^2 + \dot{y}_1^2}} + \frac{1}{2} C y_1^2\right] = 0$$
$$\frac{d}{ds}\left[(\lambda + A_{LS} + A_{SG})\frac{\dot{y}_1}{\sqrt{\dot{x}_1^2 + \dot{y}_1^2}}\right] = C y_1 \dot{x}_1$$
$$\frac{d}{ds}\left[\frac{\dot{x}_2}{\sqrt{\dot{x}_2^2 + \dot{y}_2^2}}\right] = 0 \qquad (13)$$
$$\frac{d}{ds}\left[\frac{\dot{y}_2}{\sqrt{\dot{x}_2^2 + \dot{y}_2^2}}\right] = 0.$$

*Recall that* $(x_1, y_1) = (x, y)|_{(0,l)}$ *and* $(x_2, y_2) = (x, y)|_{(l,1)}$.

The last two equations imply that $x_2(s), y_2(s)$ is a straight line. The next step is to deduce the boundary conditions and parameter $\lambda$ when the system is in equilibrium.

**Proposition 14.** *For a critical point of functional (9) the equilibrium condition at the left sheet-liquid interface (the point $(x(0), y(0))$) is that*

$$\frac{A_{LG}}{\sqrt{\dot{z}_1^2 + \dot{w}_1^2}}(\dot{w}_1, \dot{z}_1) = \frac{A_{SG} + A_{SL} + \lambda}{\sqrt{\dot{x}_1^2 + \dot{y}_1^2}}(\dot{x}_1, \dot{y}_1),$$

*where $(w_1, z_1)$ solve problem (7).*

*Proof.* Take a variation $\varphi$ of $(x_1, y_1)$ such that $\mathrm{supp}(\eta) \subset [0, l]$, Since the configuration is a minimizer, we have

$$\left.\frac{d}{dt}\right|_{t=0} E\bigg((w_1, z_1) + t\phi, (x_1, y_1) + t\eta\bigg) = 0.$$

After integrating by parts, the interior terms vanish since $(z_1, w_1, x_1, y_1)$ satisfy the Euler-Lagrange equations. We are left with the boundary term

$$(\phi_x^-, \phi_y^-) = \frac{A_{SG} + A_{SL} + \lambda}{\sqrt{\dot{x}_1^2 + \dot{y}_1^2}}(\dot{x}_1, \dot{y}_1).$$

In order to find $(\phi_x^-, \phi_y^-)$ we follow the same idea of taking a variation, integrating by parts, and keeping the boundary terms. The result is that

$$(\phi_x^-, \phi_y^-) = \frac{A_{LG}}{\sqrt{\dot{w}_1^2 + \dot{z}_1^2}}(\dot{w}_1, \dot{z}_1)$$

$\square$

Since the right hand side has norm $A_{SG} + A_{SL} + \lambda$ and the left hand side has norm $A_{LG}$, we have that $A_{SG} + A_{SL} + \lambda = A_{LG}$. The physical interpretation of this result is that the constraint of inextensibility results in a fictitious surface tension coefficient $\lambda$; and the effect of the fictitious surface tension is that the liquid-solid interface behaves just like the liquid-gas interface. We also infer that $\dot{z}_1, \dot{w}_1 = \dot{x}_1, \dot{y}_1$ at the contact point. These two observation together imply that the curve

$$\tilde{x}_1(s), \tilde{y}_1(s) = \begin{cases} x_1(l_* - s), y_1(l_* - s) \text{ if } s \in [0, l_*] \\ w_1(s - l_*), z_1(s - l_*) \text{ if } s \in [l_*, \infty) \end{cases}$$

is smooth and satisfies equation (13).

**Proposition 15.** *For a critical point of functional (5) the equilibrium condition at the contact point (the point $(x(l), y(l))$) is that*

$$\frac{A_{LG}}{\sqrt{\dot{w}_2^2 + \dot{z}_2^2}}(\dot{w}_2, \dot{z}_2) - \frac{A_{LG}}{\sqrt{\dot{x}_1^2 + \dot{y}_1^2}}(\dot{x}_1, \dot{y}_1) + \frac{2A_{SG} + \lambda}{\sqrt{\dot{x}_2^2 + \dot{y}_2^2}}(\dot{x}_2, \dot{y}_2) = 0, \qquad (14)$$

*where $(w_2, y_2)$ solves (7).*

*Proof.* The proof is similar to that of Proposition (14): by taking variations of $(x, y)$ we get that

$$(\phi_x^+, \phi_y^-) - \frac{A_{LG}}{\sqrt{\dot{x}_1^2 + \dot{y}_1^2}}(\dot{x}_1, \dot{y}_1) + \frac{2A_{SG} + \lambda}{\sqrt{\dot{x}_2^2 + \dot{y}_2^2}}(\dot{x}_2, \dot{y}_2) = 0.$$

The equation for $(\phi_x^+, \phi_y^+)$ becomes

$$(\phi_x^+, \phi_y^+) = \frac{A_{LG}}{\sqrt{\dot{w}_2^2 + \dot{z}_2^2}}(\dot{w}_2, \dot{z}_2)$$

$\square$

Finally, we recall proposition (2) to introduce a last constraint on the angles at the contact point.

**Proposition 16.** *The minimizer(s) $x(s), y(s)$ satisfy that*

$$\frac{1}{\sqrt{\dot{w}_2^2 + \dot{z}_2^2}}(\dot{w}_1, \dot{z}_1) = \frac{1}{\sqrt{\dot{x}_1^2 + \dot{y}_1^2}}(-\dot{x}_1, \dot{y}_1) \qquad (15)$$

*Proof.* Note that $\tilde{x}_1, -\tilde{y}_1$ is a minimizer of (7) such that $(\tilde{x}_1, -\tilde{y}_1) = (x_0, y_0)$ and $\lim_{s\to\infty}(\tilde{x}_1(s), \tilde{y}_1(s)) = (\infty, 0)$ by uniqueness (15) follows. $\square$

We need one last observation before we can conclude.

**Proposition 17.** *The contact point $(x_0, y_0)$ is below the critical height for the Laplace-Young equation $y^*$.*

*Proof.* Assume $y_0 > y^*$, let $l_0$ be such that $(x(l_0), y(l_0)) = (x_0, y_0)$ and let $l_*$ be such that $y(l^*) = y^*$. Let $\epsilon = l_0 - l^*$, then for $t < \epsilon$ we have that

$$\mathcal{E}(x, y, l_0) - \mathcal{E}(x, y, l_0 - t) = t(A_{LG}^* + A_{LS}^* - A_{SG}^*);$$

letting $t$ tend to 0 we obtain that $(x, y, l)$ is not a minimizer. $\square$

We finalize this section by finding the unique solution to problem (9)

**Theorem 18.** *Problem (9) has a unique solution and it is given by $x_2(s), y_2(s)$ is a vertical line, and $(\tilde{x}_1, \tilde{y}_1)$ and $(z_2, w_2)$ satisfy the (parametric) Laplace-Young equation, with height and angles given by (14) and (15).*

Even though the last result is in principle enough to find the solution of problem (9) after reparametrizing by arclength, it has an easier expression in graph coordinates $y(x)$. In this case it has an explicit solution, see for example eq 2.5 of (Anderson, Bassom, & Fowkes, 2006)

**Remark 19.** *Of course, the kink at the contact point is only an approximation. If we were to zoom in, we would observe curvature instead of a sharp angle. A natural question is then, on what scale do we observe curvature? To answer this, we write the energy functional for an interval of size $\epsilon$ around the contact point. Since the part of the sheet in contact with water meets at a non vertical angle, and since we are assuming that the derivative of the solution is smooth, we can use the implicit function theorem to work in coordinates $y(x)$, instead of $(y(s), x(s))$. Choosing a reference frame such that $x = 0$ is the contact point, we get*

$$\mathcal{E}_\epsilon = \int_{-\epsilon}^{\epsilon} h^3 E \kappa^2(x) \sqrt{1 + \left(\frac{dy}{dx}\right)^2} + \gamma \sqrt{1 + \left(\frac{dy}{dx}\right)^2} + \rho_L g y^2 dx. \qquad (16)$$

*Changing variables to $z = \epsilon x$, we get*

$$\mathcal{E}_\epsilon = \epsilon \int_{-1}^{1} \frac{h^3 E}{\epsilon^2} \kappa^2(z) \sqrt{1 + \left(\epsilon \frac{dy}{dz}\right)^2} + \gamma \sqrt{1 + \left(\epsilon \frac{dy}{dz}\right)^2} + \rho_L g y^2 dz. \qquad (17)$$

*Furthermore, we can write*

$$\begin{aligned}\rho_L g \int_{-\epsilon}^{\epsilon} y^2(x) dx &= \rho_L g \int_{-\epsilon}^{\epsilon} \left(y(0) + \int_0^x y'(\xi) d\xi\right)^2 dx \\ &= \epsilon \rho_L g y_0^2 + \mathcal{O}(\epsilon^2).\end{aligned} \qquad (18)$$

*Therefore, up to order $\epsilon$ and except for the location of the contact point, the shape of the kink is not affected by changes in $g\rho_L$. Furthermore, the bending energy becomes of the same order as other forces if $\epsilon = h^{\frac{3}{2}} \sqrt{\frac{E}{\gamma}}$.*

## Appendix

Here we prove some preliminary results used in the body of the text. We restate them for the convenience of the reader. The first is lemma 1:

**Lemma.** *There exists a $y^*$ (which depends on $A_{LG}$ and $C$) such that for every $y^* \geq y_0 > 0$ there exists a unique solution $y(x)$ to the Laplace Young equation such that $y(0) = y_0$ and $\lim_{x \to \infty} y(x) = 0$.*

*Proof.* Existence of such a function is classical; a convenient exposition is found for example in section 2a of (Anderson et al., 2006). This source also proves that solution is $C^\infty$. For uniqueness, note that solutions of the Laplace-Young equation are critical points of the functional $E_{graph}^{LY}$. Since the $(z, w) \to E_{graph}^{LY}(z, w)$ is strictly convex, the solution is unique and also a strict local minimizer. (To see that the minimizer has finite energy, notice that an admissible function is $y = \max\{y_0 - x, 0\}$, which has finite energy, therefore the minimizer has finite energy). □

Our next goal is proposition 2, which says:

**Proposition.** *The problem*

$$\min \int_0^\infty \frac{1}{2}\rho_L gy^2|\dot{x}| + A_{LG}(\sqrt{(\dot{x})^2 + (\dot{y})^2} - \dot{x})ds$$

*among $x, y$ such that $(x(0), y(0))$ is fixed and $\lim_{s \to \infty}(x(s), y(s)) = (\infty, 0)$ has a solution that is unique up to reparametrization.*

In order to prove this, we will need a technical remark about the structure of a minimizing sequence of the problem given by (7).

**Remark.** *Let $x_n(s), y_n(s)$ be a minimizing sequence. Since the functional is invariant under reparametrization, we can assume*

$$\dot{x}_n^2 + \dot{y}_n^2 = 1, \tag{19}$$

*and therefore that for $p \in (1, \infty)$ we have $x_n(s), y_n(s) \in W^{1,p}_{loc}(0, N)$. Since piecewise affine functions are dense in $W^{1,\infty}$ and the functional is l.s.c (w.r.t. strong convergence) we can assume $x_n, y_n$ is piece-wise affine. We can take a minimizing sequence that satisfies the following properties:*

1. *The function $x_n(s)$ is monotone increasing.*

2. *The function $y_n(s)$ is non negative and decreasing.*

3. *If we define $\frac{dy_n}{dx_n} = \frac{\dot{y}_n}{\dot{x}_n}$ (allowing $\dot{x}_n = 0$, in which case we get $\frac{dy_n}{dx_n} = -\infty$) then $\frac{dy_n}{dx_n}$ is monotone increasing, and because of 2 and 3, we have $\left|\frac{dy_n}{dx_n}\right|$ is decreasing.*

*Proof.* (of the remark) To establish assertion 1, suppose $x(s)$ is not monotone increasing. Then there exists $s_1 < s_2$ such that $x(s_2) < x(s_1)$. Since $x(s) \to \infty$ as $s \to \infty$ there exists $s_3 > s_2$ such that $x(s_3) = x(s_1)$, consider a trajectory $\tilde{x}(s)$ given by

$$\tilde{x}(s) = \begin{cases} (x(s), y(s)) \text{ if } s < s_1 \\ \frac{s-s_3}{s_1-s_3}(x(s_1), y(s_1))) + \frac{s_1-s}{s_1-s_3}(x(s_3), y(s_3))) \text{ if } s \in [s_1, s_3] \\ (x(s), y(s)) \text{ if } s > s_3. \end{cases}$$

Then $\tilde{x}(s)$ is increasing in the interval $(s_1, s_3)$. Since $(x(s), y(s))$ is piece-wise affine, there can be only disjoint intervals $(u_i, v_i)$ where $x(s)$ is not increasing. By applying this procedure with $s_1 = u_i$ we obtain a function $\overline{x}, \overline{y}$ such that $\overline{x}$ is monotone increasing.

Turning now to assertion 2, we argue much as we did in assertion 1. If $y(s) < 0$ for some $s$, then take $\tilde{y} = \max\{y, 0\}$. This proves we can take a minimizing sequence that is non negative. Assume there exist $s_1, s_2$ with $s_1 < s_2$ such that $0 \leq y(s_1) < y(s_2)$. Since $y(s) \to 0$ as $s \to \infty$ there exists $s_3$ such that $y(s_3) = y(s_1)$. Consider $\tilde{y}$ defined as

$$\tilde{x}(s) = \begin{cases} (x(s), y(s)) \text{ if } s < s_1 \\ \frac{s-s_3}{s_1-s_3}(x(s_1), y(s_1))) + \frac{s_1-s}{s_1-s_3}(x(s_3), y(s_3))) \text{ if } s \in [s_1, s_3] \\ (x(s), y(s)) \text{ if } s > s_3. \end{cases}$$

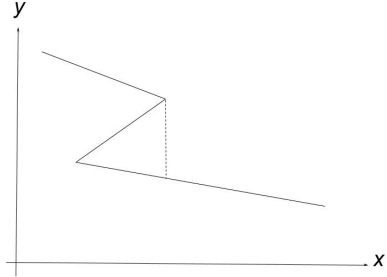

Figure 3: Rough sketch of the correction to make $x(s)$ increasing.

Then $\tilde{y}(s)$ is decreasing in the interval $(s_1, s_3)$. Since $(x(s), y(s))$ is piece wise affine, there can be only disjoint intervals $(u_i, v_i)$ where $y(s)$ is not decreasing. By applying this procedure with $s_1 = u_i$ we obtain a function $\bar{x}, \bar{y}$ such that $\bar{y}$ is monotone decreasing.

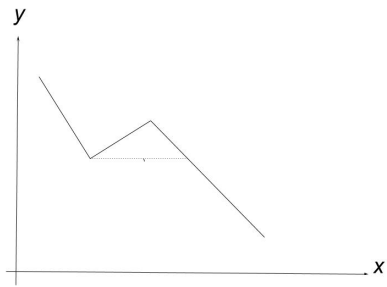

Figure 4: Rough sketch of the correction to make $y(s)$ decreasing.

Turning finally to assertion 3, we now use the hypothesis that $x(s), y(s)$ is piece-wise affine. Let $x(s), y(s)$ be given by

$$(x(s), y(s)) = \sum_i \mathbf{1}_{[a_i, a_{i+1})} \left[ (b_i s, c_i s) + (d_i, e_i) \right].$$

Then $\frac{dy}{dx} = \frac{c_i}{b_i}$ a.e. If $\frac{dy}{dx}$ is not monotone increasing, then there exists $i_0$ such that

$$\frac{c_{i_0}}{b_{i_0}} > \frac{c_{i_0+1}}{b_{i_0+1}}.$$

Define the new parameters

$$\tilde{a}_{i_0+1} = a_{i_0} - a_{i_0+2} \quad \tilde{a}_j = a_j \text{ for } j \neq i_0 + 1$$
$$\tilde{b}_{i_0} = b_{i_0+1} \quad \tilde{b}_{i_0+1} = b_{i_0} \quad \tilde{b}_j = b_j \text{ for } j \neq i_0, i_0 + 1$$
$$\tilde{c}_{i_0} = c_{i_0+1} \quad \tilde{c}_{i_0+1} = c_{i_0} \quad \tilde{c}_j = c_j \text{ for } j \neq i_0, i_0 + 1,$$

then construct the trajectory

$$\tilde{x}(s), \tilde{y}(s) = \sum_i \mathbf{1}_{[\tilde{a}_i, \tilde{a}_{i+1})} \left[ (\tilde{b}_i s, \tilde{c}_i s) + (\tilde{d}_i, \tilde{e}_i) \right], \tag{20}$$

where

$$(\tilde{d}_j, \tilde{e}_j) = (d_j, e_j) \quad \text{for } j \neq i_0, i_0 + 1$$

and $(\tilde{d}_{i_0}, \tilde{e}_{i_0}), (\tilde{d}_{i_0+1}, \tilde{e}_{i_0+1})$ are chosen so that $(\tilde{x}(s), \tilde{y}(s))$ is continuous.

Notice that the surface energy has not changed, but the gravitational energy has decreased. By applying this procedure countably many times, we obtain a function $\overline{x}, \overline{y}$ such that $\frac{d\overline{y}}{d\overline{x}}$ is monotone increasing. The proof of the remark 1 is now complete.

$\square$

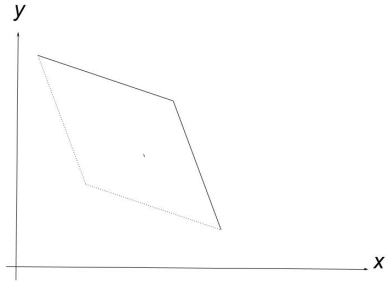

Figure 5: Rough sketch of the correction to make $\frac{dy}{dx}(s)$ increasing.

We now use the preceding remark to prove Proposition 2.

*Proof.* (of Proposition 2) Because of (19) we have that $x_n, y_n$ is bounded in $W^{1,\infty}[0, N]$ for every $N$, hence we can extract a sequence of sequence $x_{n_{j_N}}, y_{n_{j_N}}$ such that $(x_{n_{j_N}}, y_{n_{j_N}})$ converges weakly in $W^{1,p}[0, N]$ and $x_{n_{j_{N+1}}}, y_{n_{j_{N+1}}}$ is a subsequence of $x_{n_{j_N}}, y_{n_{j_N}}$. Hence we can extract a diagonal subsequence that converges weakly in $W^{1,p}[0, N]$ for all $N$. Let $(x(s), y(s))$ be the limit. Since the functional is invariant under reparametrization, we can take $x(s), y(s)$ parametrized by arclength.

For each $n$ we have that $x_n(s) \geq s - y_0$, hence $x(s) \geq s - y_0$ and $x(s) \to \infty$. Since $\dot{x}_n$ is increasing and positive, we have that $\dot{x}$ is non decreasing and non negative. Since $x(s) \to \infty$ as $s \to \infty$ then $\dot{x}$ can only be 0 in an interval $(0, \epsilon]$. Therefore in the interval $(\epsilon, \infty)$ we have that $y = y(x)$ (strictly speaking, there exists $f$ such that for $s > \epsilon$ we have $y(s) = f(x(s))$.) Using the previous results, we have that $\frac{dy}{dx}$ is decreasing and non positive. Since $(x(s), y(s))$ is a critical point of (2), we have that $y(x)$ is a critical point (and therefore minimizer) of (7), and therefore $y(x)$ solves

$$\frac{y''}{\sqrt{(1 + y'^2)^3}} = y. \tag{21}$$

Since $y(x)$ and $y'(x)$ are monotone, we have that $y''(x) \to 0$ as $x \to \infty$. Therefore $y(x) \to 0$ as $x \to \infty$.

We now split into two cases. In the case $y_0 < y^*$ we show that $\lim_{x \to 0^+} = y_0$, which is equivalent to showing that $\epsilon = 0$. Note that $\lim_{x \to 0^+} \neq \infty$ because $y_0 < y^*$. Assume that $\epsilon \neq 0$, the idea to get a contradiction is to chop off a sequence of triangles that decrease the energy. For a given $\delta < 0$ define the points

$$\mathbf{x}_{\delta+} = (x(\epsilon + \delta), y(\epsilon + \delta)) \ \mathbf{x}_{\delta-} = (x(\epsilon - \delta), y(\epsilon - \delta));$$

define the function $x_\delta(s), y_\delta(s)$ as

$$x_\delta(s), y_\delta(s) = \begin{cases} x(s), y(s) \text{ if } s \in [0, \epsilon - \delta) \\ (-\frac{s - [\epsilon + \delta]}{2\delta})\mathbf{x}_{\delta-} + (\frac{s - (\epsilon - \delta)}{2\delta})\mathbf{x}_{\delta+} \text{ if } s \in [\epsilon - \delta, \epsilon + \delta] \\ x(s), y(s) \text{ if } s > \epsilon + \delta. \end{cases}$$

Then the increase in gravitational energy (the $y^2 \dot{x}$ term) is $\mathcal{O}(\delta^2)$, but the decrease in surface energy (the $\sqrt{\dot{x}^2 + \dot{y}^2} - \dot{x}$ term) is $\mathcal{O}(\delta)$, by taking $\delta$ small enough we get that $x(s), y(s)$ is not a critical point of the functional.

In the case $y_0 > y^*$, we have that $\epsilon = y_0 - y^*$. To prove this, we proceed by contradiction, if $y_0 - y^* < \epsilon$, then by the same reasoning, $x(s), y(s)$ is a graph for $s > \epsilon$ and then there is a corner at $x(\epsilon), y(\epsilon)$. We can proceed similarly to case 1 and obtain that $x(s), y(s)$ is not a critical point. If $y_0 - y^* > \epsilon$ then we can write with an abuse of notation $y = y(x)$ for $s > \epsilon$. Then $y$ solves (21), and $y(0) > y^*$, this is a contradiction. □

Next, we provide the proof of Proposition 3, it states:

**Proposition.** *The function $\phi^\pm$ (defined by (6)) is smooth away from $y^*$.*

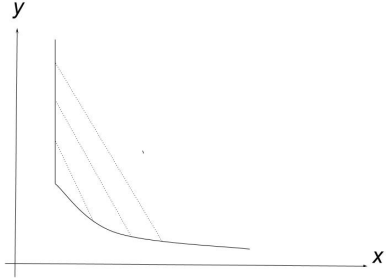

Figure 6: Rough sketch of the argument.

*Proof.* Let $y(x)$ and $\tilde{y}(x)$ be solutions such that $y(0) = y_0$ and $\tilde{y}(0) = \tilde{y}_0$ with $\lim_{x\to\infty} y(x) = \lim_{x\to\infty} \tilde{y}(x) = 0$. Assume WLOG that $\tilde{y}_0 > y_0$. By continuity there exists $\tau$ such that $\tilde{y}(\tau) = y_0$. By uniqueness of the solution it follows that $y(x) = \tilde{y}(x + \tau)$.) Let $Y(x)$ be the unique solution to the LY equation such that $Y(0) = y^*$. Then for every $0 < \rho < y^*$ the solution to the LY with decay at infinity and initial condition $y(0) = \rho$ can be written as

$$y(x) = Y(x + Y^{-1}(\rho)).$$

Hence

$$\phi^\pm(y_0) = \int_0^\infty E^{LY}_{graph}(y, y')dx$$
$$= \int_{Y^{-1}(y_0)}^\infty E^{LY}_{graph}(Y, Y')dx.$$

We have proved that $\phi^\pm$ is continuous. Furthermore, the smoothness of $Y(x)$ together with IFT imply that $\phi^\pm$ is smooth away from $y = y*$. □

Finally, we provide the proof of lemma 7. It states:

**Lemma.** *Let $f : [0, 1] \to \mathbf{R}^2$ be a $C^\infty$ short embedding. For every $\epsilon > 0$ there exists $g \in C^\infty[0, 1]$ such that*

- *The function $g$ is an isometric embedding.*
- *$|f(s) - g(s)| < \epsilon$ for all $s \in [0, 1]$.*

*Proof.* We shall proceed in 3 distinct steps.

*Step 1* Let $f : [0, 1] \to \mathbf{R}^2$ be a short embedding. Then there exists an absolute constant $K$ such that for any isometric immersion $g : [0, 1] \to \mathbf{R}^2$ such that $g(0) = f(0)$, $g(1) = f(1)$ and we have

$$\| f - g \|_{C^0} \leq K. \tag{22}$$

Also, there exists an isometric immersion $g : [0, 1] \to \mathbf{R}^2$ such that

$$g^{(n)}(0) - f^{(n)}(0) = g^{(n)}(1) - f^{(n)}(1) = 0 \text{ for } n = 0, 1, 2... \tag{23}$$

To establish (22) we simply observe that for any such isometric embedding $g$ we have

$$\| g(x) - f(x) \| \leq \| g(x) - g(0) \| + \| f(0) - f(x) \|$$
$$\leq 2.$$

For (23), define $\tilde{f} : [0, l] \to \mathbf{R}^2$ the reparametrization by arclength of $f$, where $l$ is the length of the curve. Let $\phi : [0, l] \to \mathbf{R}^2$ be defined as

$$\phi(x) = \frac{(\dot{\tilde{f}}(x))^\perp}{\| (\dot{\tilde{f}}(x))^\perp \|}.$$

Let $\mu : [0, l] \to \mathbf{R}$ be such that

$$\mu(x) \geq 0$$
$$\mu(x) > 0 \text{ in } (0, l)$$
$$\mu^{(n)}(0) = \mu^{(n)}(l) = 0 \text{ for } n = 0, 1, 2...$$

Define the function
$$\tilde{g}(x) = \tilde{f}(x) + C\mu(x)\phi(x), \tag{24}$$

where $C$ is such that

$$\int_0^l \| \frac{d}{dx}\left(\tilde{f}(x) + C^2\mu(x)\phi(x)\right) \| \, dx = 1.$$

Note that $C$ exists by Intermediate Value Theorem. Let $g$ be a reparametrization of $\tilde{g}$ by arclength, then $g$ is the function we are looking for.

*Step 2* We claim that if $f : [0, \tau] \to \mathbf{R}^2$ is a short embedding, then there exists an absolute constant $K$ such that for any isometric immersion $g : [0, \tau] \to \mathbf{R}^2$ such that $g(0) = f(0)$ and $g(\tau) = f(\tau)$ we have

$$\| f - g \|_{C^0} \leq K\tau;$$

also, there exists an isometric immersion $g : [0, \tau] \to \mathbf{R}^2$ such that

$$g^{(n)}(0) - f^{(n)}(0) = g^{(n)}(1) - f^{(n)}(1) = 0 \text{ for } n = 0, 1, 2...$$

The proof is easy: just apply Step 1 one to $\tilde{f} : [0, 1] \to \mathbf{R}^2$ defined as

$$\tilde{f}(x) = \frac{1}{\tau}f(\tau x)$$

*Step 3* We now prove the assertion of the Proposition. Let $N$ be such that $\frac{K}{N} < \epsilon$. Let $f_j = f|_{[\frac{j-1}{n}, \frac{j}{n}]}$. By step 2 there exists $g_j : [\frac{j-1}{n}, \frac{j}{n}] \to \mathbf{R}^2$ such that $g_j$ is an isometric embedding and

$$g_j^{(n)}(\frac{j-1}{n}) - f_j^{(n)}(\frac{j-1}{n}) = g_j^{(n)}(\frac{j}{n}) - f_j^{(n)}(\frac{j}{n}) = 0 \text{ for } n = 0, 1, 2...$$

Notice that the immersion $g_j$ is am embedding for $\epsilon$ small enough by a direct application of the tubular neighborhood theorem. We also have

$$\| f_j - g_j \|_{C^0} \leq \frac{K}{n} < \epsilon.$$

Therefore $g$ defined as

$$g(x) = g_j(x) \text{ if } x \in [\frac{j-1}{n}, \frac{j}{n}]$$

is the function we are looking for. □

## Acknowledgments

The author wants to thank Benny Davidovitch for introducing him to this topic.